\documentclass[12pt,twoside]{article}		
\input{atmp.sty}
\preptrue
\copyrightnotice{2001}{5}{1000}{1012}

\begin{document}


\numberwithin{equation}{section}

\newcommand{\nc}{\newcommand}
\nc{\ona}{\operatorname}
\nc{\al}{\alpha}
\nc{\De}{\Delta}
\nc{\e}{\varepsilon}
\nc{\ga}{\gamma}
\nc{\Ga}{\Gamma}
\nc{\la}{\lambda}
\nc{\La}{\Lambda}
\nc{\Om}{\Omega}
\nc{\si}{\sigma}
\nc{\Si}{\Sigma}
\nc{\va}{\varphi}
\nc{\f}{\frac}
\nc{\iy}{\infty}
\nc{\pa}{\partial}
\nc{\na}{\nabla}
\nc{\ra}{\rightarrow}
\nc{\wt}{\widetilde}
\nc{\Area}{{\rm Area}\,}
\nc{\Ric}{{\rm Ric}\,}

\def\<{{\langle}}
\def\>{{\rangle}}

\newtheorem{theorem}{Theorem}[section]
\newtheorem*{remark}{Remark}

\title{Geometry of three manifolds and existence of Black Hole due to boundary effect}
\author{Shing Tung Yau}

\address{Department of Mathematics\\
Harvard University\\
Cambridge, MA 02138}
\addressemail{yau@math.harvard.edu}

\url{xxxx}

\section{Introduction}

In this paper, we observe that the brane functional studied in
\cite{5} can be used to generalize some of the works that Schoen and I
\cite{4} did many years ago. The key idea is that if a three dimensional
manifold $M$ has a boundary with strongly positive mean curvature, the
effect of this mean curvature can influence the internal geometry of
$M$. For example, if the scalar curvature of $M$ is greater than
certain constant related to this boundary effect, no
incompressible surface of higher genus can exist.
\newpage

A remarkable statement in general relativity is that if the mean
curvature of $\pa M$ is strictly greater than the trace of $p_{ij}$
(the second fundamental form of $M$ in space time), the value of this
difference can provide the existence of apparent horizon in $M$.
In fact, matter density can even be allowed to be negative if this
boundary effect is very strong. Theorem 5.2 is the major result of
this paper.

\pagestyle{myheadings}
\markboth{\hfil{\sc \small S.-T. YAU\hspace*{30mm}}\hfil}%
{\hfil{\sc \small\hspace*{10mm} GEOMETRY of THREE MANIFOLDS}\hfil}

\section{Existence of stable incompressible surfaces with constant
mean curvature}

We shall generalize some of the results of Schoen-Yau \cite{1} and
Meeks-Simon-Yau \cite{3}.

Let $M$ be a compact three dimensional manifold whose boundary $\pa M$
has mean curvature (with respect to the outward normal) not less than
$c>0$. Assume the volume form of $M$ can be written as $d\La$ where
$\La$ is a smooth two form.

Let $f:\Si\ra M$ be a smooth map from a surface $\Si$ into $M$ which
is one to one on $\Pi_1(\Si)$. We are interesting in minimizing the
energy
$$
E_c(f)=\f{1}{2}\int_{\Si}|\na f|^2-c\int_{\Si}f^*\La.
$$

There are two different hypothesis we shall make for the existence of
surfaces which minimizes $E_c$.

\begin{theorem}			
Assume the existence of embedded $\Si$ with $\Pi_1(\Si)\ra\Pi_1(M)$ to
be one to one. Assume that for any ball $B$ in $M$, 
the volume of the ball $B$ is not greater than $\f{1}{c}\,\Area(\pa
B)$. Then we can find a surface isotopic to $\Si$ which minimizethe
functional $\Area(\Si)-c\int_{\Si}\La$.
\end{theorem}

\begin{proof}
This follows from the argument of Meeks-Simon-Yau \cite{3}. The
hypothesis is method to deal with the cut and paste argument.
\end{proof}

\begin{theorem}			
Assume that the supremum norm of $\La$ is not greater than
$c^{-1}$. Then for any $c'<c$, we can find a conformal map from some
conformal structure on $\Si$ to $M$ which induces the same map is 
$f_*=\Pi_1(\Si)\ra \Pi_1(M)$ and minimize the energy
$\f{1}{2}\int_{\Si}|\na f|^2-c'\int_{\Si}f^*\La$.
\end{theorem}

\begin{proof}
Since $c'<c$, the energy is greater than a positive multiple of the
standard energy. Hence the argument of Schoen-Yau \cite{1} works.
\end{proof}

\begin{remark}
{\em
It should be possible to choose $c'=c$ in this last theorem.
}
\end{remark}

\section{Second variational formula}

Let $\Si$ be the stable surface established in section 2. Then the
variational formula shows that the mean curvature of $\Si$ is equal to
$c$. The second variational formula shows that for all $\va$ defined
on $\Si$,
\begin{equation}
\int_{\Si}|\na\va|^2-\int_{\Si}\left(\Ric_M(\nu,\nu)
+\Si h^2_{ij}\right)\va^2\geq 0
\end{equation}
where $\Ric_M(\nu,\nu)$ is the Ricci curvature of $M$ along the normal
of $\Si$ and $h_{ij}$ is the second fundamental form.

The Gauss equation shows that
\begin{equation}
\Ric_M(\nu,\nu)=\f{1}{2}R_M-K_{\Si}+\f{1}{2}(H^2-\Si h^2_{ij})
\end{equation}
where $R_M$ is the scalar curvature, $K_{\Si}$ is the Gauss curvature
of $\Si$ and $H=$ trace of $h_{ij}$ is the mean curvature.

Hence
\begin{equation}	
\int_{\Si}|\na\va|^2\geq \f{1}{2}\int_{\Si}(R_M+\Si h_{ij}^2+H^2)\va^2
-\int_{\Si}K_{\Si}\va^2.
\end{equation}

Since $\Si h^2_{ij}\geq \f{1}{2}H^2$, we conclude that
\begin{equation}	
\int_{\Si}|\na\va|^2
\geq\f{1}{2}\int_M\left(R_M+\f{3}{2}H^2\right)\va^2
-\int_{\Si}K_{\Si}\va^2.
\end{equation}

If $\chi(\Si)\leq 0$, we conclude by choosing $\va =1$, that
\begin{equation}	
\int_{\Si}\left(R_M+\f{3}{2}c^2\right)\leq 0.
\end{equation}

\begin{theorem}	
If $R_M+\f{3}{2}c^2\geq 0$, any stable orientable surface $\Si$ with
$\chi(\Si)\leq 0$ must have $\chi(\Si)=0$ and $R_M+\f{3}{2}c^2=0$ along
$\Si$. Furthermore $\Si$ must be umbilical.
\end{theorem}

Let us now see whether stable orientable $\Si$ with
$\chi(\Si)=0$ can exist or not.
If $R_M+\f{3}{2}c^2>0$ at some point of $M$ and $R_M+\f{3}{2}c^2 \geq 0$
everywhere. We can deform the metric conformally so that $R_M
+\f{3}{2}c^2>0$ everywhere while keeping mean curvature of $\pa M$ not
less than $c$. (This can be done by arguments of Yamabe problem.)
In this case, incompressible torus does not exist.

Hence we may assume $R_M+\f{3}{2}c^2=0$ everywhere. In this case, we
deform the metric to $g_{ij}-t(R_{ij}-\f{R}{3}g_{ij})$. By
computation, one sees that unless $R_{ij}-\f{R}{3}g_{ij}$ everywhere,
the (new) scalar curvature will be increased.

Let $\Si$ be the stable surface with constant mean curvature with
respect to the metric $g_{ij}$. We can deform the surface $\Si$ along
the normal by multiplying the normal with a function $f$. For this
surface $\Si_f$, we look at the equation $H_{t}(\Si_f)=c$ where
$H_{t}$ is the mean curvature with respect to the new metric at time
$t$. As a function of $t$ and $f$, $H_t(\Si_f)$ define a mapping
into the Hilbert space of functions on $\Si$. The linearized operator
with respect to the second (function) variable is
$-\Delta-(\Ric(\nu,\nu)+\Si h_{ij}^2)$. This operator is self-adjoint
and if there is no kernel, we can solve the equation
$H_t(\Si_f)=c$ for $t$ small.

We conclude that if $-\Delta-(\Ric(\nu,\nu)+\Si h_{ij}^2)$ has no
kernel and if the metric is not Einstein, we can keep mean curvature
constant and scalar curvature greater than $-\f{3c^2}{2}$. On the
other hand, if the metric is Einstein, we can use argument in \cite{5}
to prove that $M$ is the warped product of the flat torus with $R$.

If $-\Delta-(\Ric(\nu,\nu)+\Si h_{ij}^2)$ has kernel, it must be a
positive function $f$ defined on $\Si$. (This comes from the fact that it
must be the first eigenfunction of the operator.) Hence
\begin{align*}
\Delta(\log f)+|\na \log f|^2
&=-(\Ric(\nu,\nu)+\Si h_{ij}^2)\\
&\leq -\f{1}{2}\left(R_M+\f{3}{2}H^2\right)+K_{\Si}.
\end{align*}

Since
$$
\int_{\Si}\left(R_M+\f{3}{2}H^2\right)\geq 0
$$
and
$$
\int_{\Si}K_{\Si}=0,
$$
$f$ must be a constant and 
$$
K_{\Si}\geq -\f{1}{2}\left(R_M+\f{3}{2}H^2\right)= 0.
$$

Hence $K_{\Si}=0$ and $R_M=-\f{3}{2}H^2$ is constant along $\Si$. Also
$h_{ij}=\f{H}{2}g_{ij}$ and $\Ric(\nu,\nu)+\Si h_{ij}^2=0$ along
$\Si$. 

If we compute the first order deformation of the mean curvature of
$\Si$ along the normal, it is trivial as
$h_{ij}=\f{H}{2}g_{ij}$, $R=-\f{3}{2}H^2$ and $R(\nu,\nu)=-\Si
h_{ij}$.

In conclusion, the mean curvature is equal to $H$ up to first order in
$t$ while we can increase the scalar curvature of $M$ up to first
order (unless $R_{ij}=\f{R}{3}g_{ij}$ everywhere). We can therefore
prove the following

\begin{theorem}	
Let $M$ be a three dimensional complete manifold with scalar curvature
not less than $-\f{3}{2}c^2$ and one of the component of $\pa M$ is an
orientable incompressible surface with nonpositive Euler number and
mean curvature $\geq c$. Suppose that for any ball $B$ in $M$, the
area of $\pa B$ is not less than $c{\rm Vol}\,(B)$. Then $M$ is
isometric to the warped product of the flat torus with a half line.
\end{theorem}

\section{Geometry of manifolds with lower\\ bound on scalar curvature.}

In this section, we generalize the results of Schoen-Yau \cite{4}.

Given a region $\Om$ and a Jordan curve $\Ga\subset \pa\Om$ which
bounds an embedded disk in $\Om$ and a subdomain in $\pa\Om$ we define
$R_{\Ga}$ to be the supremium of $r>0$ so that $\Ga$ does not bound a disk
inside the tube of $\Ga$ with radius $r<R_{\Ga}$. We define
Rad$\,(\Om)$ to be the supremium of all such $R_{\Ga}$. (Note that this
concept can be generalized to higher homology or homotopic groups. The
Radius that is defined in this manner will be sensitive to geometry of
higher homology or higher homotopic groups.)

Let $R$ be the scalar curvature of $M$ and $h$ be a function defined
on $M$ and $k$ be a function defined on $\pa M$ so that for any smooth
function $\va$
\begin{equation}
\int_M|\na\va|^2+\f{1}{2}\int_MR\va^2+\int_{\pa M}k\va^2
\geq \int_Mh\va^2.
\end{equation}

Let $f$ be the positive first eigenfunction of the operator
$-\Delta+\f{1}{2}R-h$ so that
\begin{equation}
\begin{cases}
-\Delta f+\f{1}{2}Rf-hf=\la f&\\
\f{\pa f}{\pa\nu}+kf=0&\quad \text{on }\;\pa M.
\end{cases}
\end{equation}

Let $\Ga$ be a Jordan curve on $\pa\Om$ which defines Rad$\,(\Om)$
up to a small constant.
Let $\Si$ be a disk in $\Om$ with boundary $\Ga$ such that $\Si$
together with a region on $\pa\Om$ bounds a region $\Om_{\Si}$.

Assume that $\pa\Om$ has mean curvature $H$ so that $f(H-k)$ is
greater than $cf$. Then we define a functional
\begin{equation}
L(\Si)=\int_{\Si}f-c\int_{\Om_{\Si}}f.
\end{equation}

Let us now demonstrate that $\pa\Om$ forms a ``barrier''for the
existence of minimum of $L(\Si)$.

Let $r$ be the distance function to $\pa\Om$. Let us assume that $\Ga$
is in the interior of $\Om$. If $\Si$ touches $\pa\Om$, we look at the
domain $\Om_{\Si}\cap\{0<r<\e\}=\Om_{\Si,\e}$. Then
\begin{equation}
\int_{\pa(\Om_{\Si,\e})}f\f{\pa r}{\pa\nu}
=\int_{\Om_{\Si,\e}}f\Delta r+\int_{\Om_{\Si,\e}}\na f\cdot\na r.
\end{equation}

When $\e$ is small, $f\Delta r+\na f\cdot \na r$ is close to the
boundary value $-\f{\pa f}{\pa\nu}-Hf$ on $\pa\Om$. Hence
\begin{equation}
\int_{\pa(\Om_{\Si,\e})}f\f{\pa r}{\pa\nu}
<-\int_{\Om_{\Si,\e}}cf.
\end{equation}

Since $\left|\f{\pa r}{\pa\nu}\right|\leq 1$ and $\f{\pa r}{\pa\nu}=1$
along $r=\e$, we conclude that if we replace $\Si$ by
$(\Si\smallsetminus \pa\Om_{\Si,\e})\cup
(\pa\Om_{\Si,\e}\cap\{r=\e\})$, then the new surface will have
strictly less energy than $L_f(\Si)$. Hence when we minimize $L$,
$\pa\Om$ forms a barrier.

By standard geometric measure theory, we can find a surface $\Si$
which minimize the functional $L_f$. (We start to minimize the
functional $\int_{\Si}f-tc\int_{\Om_{\Si}}f$ when $t$ is small.)

For this surface, we can compute both the first variational and second
variational formula and obtain from the first variational formula
\begin{equation}
\f{\pa f}{\pa\nu}+Hf=cf.
\end{equation}

The second variational formula has contributions from two terms. The
second term gives rise to
\begin{equation}
\int c\left(\f{\pa f}{\pa\nu}+Hf\right)=c^2\int_{\Si}f.
\end{equation}

Using (4.6) the first term of the second variational formula gives
\begin{align}
0&\leq \int_{\Si}|\na\va|^2f\\
&\quad-\int_{\Si}\left(\f{1}{2}R_M-K_{\Si}\right)\va^2f
-\int\det(h_{ij})\va^2f\nonumber\\
&\quad+\int_{\Si}\left(\Delta_Mf-\Delta_{\Si}f
-H\f{\pa f}{\pa\nu}\right)\va^2
-\int(\Si h^2_{ij}-H^2)\va\nonumber\\
&\quad+2\int_{\Si}\f{\pa f}{\pa\nu} H\va^2-\int_{\Si}c^2f\va^2\nonumber\\
&\leq \int_{\Si}|\na\va|^2
+\int_{\Si}\left(\Delta_Mf-\f{1}{2}R_Mf\right)\va^2\nonumber\\
&\quad-\int_{\Si}(\Delta_{\Si}f-K_{\Si}f)\va^2
+\f{1}{4}\int_{\Si}H^2f\va^2
+\int_{\Si}\f{\pa f}{\pa\nu}H\va^2
-c^2\int_{\Si}f\va^2\nonumber
\end{align}
where $\Delta_M$ and $\Delta_{\Si}$ are the Laplacian of $M$ and $\Si$
respectively and $\va$ is any function vanishing on $\pa\Si$.

We conclude from $\f{\pa f}{\pa\nu}+Hf=cf$ that
\begin{align}
&\f{1}{4}\int_{\Si}H^2f\va^2+\int_{\Si}\f{\pa f}{\pa\nu}H\va^2
-c^2\int_{\Si}f\va^2\\
&\leq -\f{3}{4}\int_{\Si}H^2f\va^2+c\int_{\Si}Hf\va^2
-c^2\int_{\Si}f\va^2\nonumber\\
&\leq -\f{2}{3}c^2\int_{\Si}f\va^2,\nonumber
\end{align}
\begin{equation}
\int_{\Si}|\na\va|^2f-\int_{\Si}(\De_{\Si}f-K_{\Si}f)\va^2
-\int_{\Si}\left(h+\la+\f{2c^2}{3}\right)f\va^2\geq 0
\end{equation}
where $\la$ is the first eigenvalue of the operator $-\De+\f{R}{2}-h$
with boundary value given by $\f{\pa f}{\pa\nu}+kf=0$.

By the argument of \cite{4}, we see that for any point $p\in\Si$,
there exists a curve $\si$ from $p$ to $\pa\Si$ with length $l$ such
that 
\begin{equation}
\int_0^l\left(h+\la+\f{2c^2}{3}\right)\va^2
\leq\f{3}{2}\int_0^l(\va')^2
\end{equation}
where $l$ is the length of the curve $\si$ and $\va$ vanishes at 0 and
$l$.

\begin{theorem}		
Let $M$ be a three dimensional manifold so that (4.1) holds. Let
$\la$ be the first eigenvalue of the operator (4.2). Suppose that
the mean curvature of $\pa M$ minus $k$ is greater than a constant
$c>0$. Then for any closed curve $\Ga\subset M$, there is a surface
$\Si$ that $\Ga$ bounds in $M$ so that for any point $p\in\Si$, there
is a curve $\si$ from $p$ to $\pa\Si$, inequality (4.11) holds.
\end{theorem}

\section{Existence of Black Holes}

For a general initial data set for the Einstein equation, we have two
tensors $g_{ij}$ and $p_{ij}$. The local energy density and linear
momentum are given by
\begin{align}
\mu&=\f{1}{2}\left[R-\Si p^{ij}p_{ij}+\left(\Si p_i^i\right)\right]\\
J^i&=\sum_jD_j\left[p^{ij}-\left(\Si p^k_k\right)
g^{ij}\right]\nonumber
\end{align}

In \cite{2}, Schoen and I studied extensively the following equation
initiated by Jung
\begin{equation}
\sum_{i,j}\left(g^{ij}-\f{f^if^j}{1+|\na f|^2}\right)
\left(\f{D_iD_jf}{(1+|\na f|^2)^{1/2}}-p_{ij}\right)=0.
\end{equation}

For the metric
$$
\wt{g}_{ij}=g_{ij}+\f{\pa f}{\pa x^i}\;\f{\pa f}{\pa x^j},
$$
one has the following inequality
\begin{align}
2(\mu-|J|) &\leq \bar{R}-\sum_{i,j}(h_{ij}-p_{ij})^2\\
&\quad-2\sum(h_{i4}-p_{i4})^2+2\sum D_i(h_{i4}-p_i).\nonumber
\end{align}

Hence for any function $\va$
\begin{align}
2\int_M(\mu-|J|)\va^2
&\leq \int_M\bar{R}\va^2
-2\int_M\sum(h_{i4}-p_{i4})^2\\
&\quad-\int_M4\va(\na_i\va)(h_{i4}-p_{i4})
+2\int_{\pa M}(h_{\nu4}-p_{\nu4})\va^2\nonumber\\
&\leq \int_M\bar{R}\va^2+2\int_M|\na\va|^2
+2\int_{\pa M}(h_{\nu4}-p_{\nu4})\va^2.\nonumber
\end{align}

Hence in (4.1) we can take
\begin{equation}
h=(\mu-|J|),\qquad k=h_{\nu4}-p_{\nu4}.
\end{equation}

Let $e_1,e_2,e_3$ and $e_4$ be orthonormal frame of the graph so that
$e_1,e_2$ is tangential to $\pa M$ (assuming $f=0$ on $\pa M$) and
$e_3$ is tangential to the graph but normal to $\pa M$. By assumption,
\begin{equation}
h_{\nu4}=h_{34}=\<\na_{e_4}e_4,e_3\>.
\end{equation}

Let $w$ be the outward normal vector of $\pa M$ in the horizontal space
where $f=0$. Hence
\begin{align}
h_{34} &=\<e_4,w\>\<\na_we_4,e_3\>\\
&=-\<e_4,w\>\<e_4,\na_we_3\>\nonumber\\
&=\f{-\<e_4,w\>}{\<e_3,w\>}\;\<e_4,\na_{e_3}e_3\>.\nonumber
\end{align}

Since $-\sum^3_{i=1}\<e_4,\na_{e_i}e_i\>$ is the mean curvature of the
graph of $f$ which is tr$\,p$, we conclude that
\begin{align}
h_{34} &=\f{\<e_4,w\>}{\<e_3,w\>}\;
\left({\rm tr}\,p+\sum_{i=1}^2\<e_4,\na_{e_i}e_i\>\right)\\
&=\f{\<e_4,w\>}{\<e_3,w\>}\;
\left({\rm tr}\,p+\<e_4,w\>\sum_{i=1}^2\<w,\na_{e_i}e_i\>\right).\nonumber
\end{align}

The mean curvature of $\pa\Om$ with respect to the metric $g_{ij}+
\f{\pa f}{\pa x_i}\,\f{\pa f}{\pa x_j}$ is given by
$$
-\sum_{i=1}^2\<\na_{e_i}e_i,e_3\>
=-\sum_{i=1}^2\<\na_{e_i}e_i,w\>\,\<w,e_3\>.
$$

Hence the difference between mean curvature and $k$ is given by
\begin{equation}
-\f{\<e_4,w\>}{e_3,w\>}\,{\rm tr}\,p
+p_{34}+\left(\<w,e_3\>+\f{\<w,e_4\>^2}{\<w,e_3\>}\right)H_{\pa\Om}
\end{equation}
where $H_{\pa\Om}$ is the mean curvature of $\pa\Om$ with respect to
the metric $g_{ij}$.

Since 
\begin{align*}
p_{34} &=\<e_4,w\>\;p(e_3,w)\\
&=\f{\<e_4,w\>}{\<e_3,w\>}\;p\<e_3,e_3\>.
\end{align*}

We conclude that 
the expression (5.9) is given by
\begin{align}
&-\f{\<e_4,w\>}{\<e_3,w\>}\,({\rm tr}\,_{\pa\Om}p)
+\f{1}{\<e_3,w\>}H_{\pa\Om}\\
&\geq \left(H_{\pa\Om}-|{\rm tr}\,_{\pa\Om}p|\right)
\<e_3,w\>^{-1}.\nonumber
\end{align}

We shall assume $H_{\pa\Om}>|{\rm tr}\,_{\pa\Om}p|$ and we can choose
$c$ to be lower bound of $H_{\pa\Om}-|{\rm tr}\,_{\pa\Om}p|$.

We need to solve the Dirichlet problem for $f$ with $f=0$ on
$\pa\Om$. While most of the estimates were made in \cite{2}, we need
to construct a barrier for the boundary valued problem.

Let $\va$ be an increasing function defined on the interval $[0,\e]$
so that $\va'(\e)=\iy$. Let $d$ be the distance function from $\pa\Om$
measured with respect to $g_{ij}$.

Then $\va(d)$ can be put in (5.2) and when $\e$ is small, we obtain
the expression
\begin{equation}
\f{\va'}{\sqrt{1+(\va')^2}}\,\left(-H_{\pa\Om}\right)
-{\rm tr}_{\pa\Om}p
+\f{\va''}{(1+\va^{\prime 2})^{3/2}}-\f{p_{\nu\nu}}{1+\va^{\prime 2}}.
\end{equation}

To construct a supersolution, we need this expression to be
nonpositive. When $\e$ is small, and $\va'$ is very large, the
condition is simply $H_{\pa\Om}>{\rm tr}_{\pa\Om}p$. Similarly, we can
construct a subsolution using $-\va(d)$. The conclusion is that we can
solve (5.2) if $H_{\pa\Om}>|{\rm tr}_{\pa\Om}p|$. We have therefore
arrived at the following conclusion

\begin{theorem}		
Let $M$ be a space like hypersurface in a four dimensional
spacetime. Let $g_{ij}$ be the induced metric and $p_{ij}$ be the
second fundamental form. Let $\mu$ and $J$ be the energy density and
local linear momentum of $M$. Suppose the mean curvature $H$ of $\pa
M$ is greater than tr$\,_{\pa M}(p)$. Assume that $H-|{\rm tr}_{\pa M}p|
\geq c\geq 0.$ Let $\Ga$ be a Jordan curvature in $\pa M$ that bounds
a domain in $\pa M$. If $M$ admits no apparent horizon, then there
exists a surface $\Si$ in $M$ bounds by $\Ga$ so that for any point
$p\in\Si$, there is a curve $\si$ with length $l$ from $p$ to $\Ga$
and
\begin{equation}
\int_0^l\left((\mu-|J|)+\f{3}{2}c^2\right)\va^2ds
\leq\int_0^l|\na\va|^2ds
\end{equation}
where $\va$ is any function vanishing at $0$ and $l$.
\end{theorem}

\begin{theorem}		
Let $M$ be a space like hypersurface in a spaetime. Let $g_{ij}$ be its
induced metric and $p_{ij}$ be its second fundamental form. Assume
that the mean curvature $H$ of $\pa M$ is strictly greater than $|{\rm
tr}_{\pa m}(p)|$. Let $c=\min(H-|{\rm tr}_{\pa m}(p)|)$ if
${\rm Rad}\,(M)\geq \sqrt{\f{3}{2}}\,\f{\Pi}{\sqrt{\La}}$ where $\La=
\f{2}{3}c^2+\mu-|J|$, then  $M$ must admit apparent horizons in its interior.
\end{theorem}

An important point here is that the curvature 
$H-|{\rm tr}_{\pa m}(p)|)$ of the boundary itself can give rise to
Black Hole.

The inequality actually shows that as long as $\mu-|J|\geq 0$
everywhere, $\f{2c^2}{3}+\mu-|J|$ to be large in a reasonable ringed
region and Rad$\,(\Om)$ is large, an apparent horizon will form in $M$.

\end{document}